\documentclass{article}
\usepackage{CJKutf8}

\usepackage{amsmath}
\usepackage{caption}
\usepackage{graphicx}
\usepackage{subfigure}
\usepackage{wrapfig}
\usepackage{float}

\usepackage{picinpar}

\begin{document}
\begin{CJK*}{UTF8}{gbsn}

\renewcommand{\baselinestretch}{1.2}
\voffset=-1cm \hoffset=0cm \textwidth=18cm \textheight=22cm
\topmargin=0cm
\parskip=5pt

\makeatletter 
\renewcommand{\@thesubfigure}{\hskip\subfiglabelskip}
 \makeatother

\title{The Beginning and the Unfinished Story of the Millennium Prize Problems\\  \large 千禧年大奖难题之始与未终}

\author{Arthur Jaffe\footnote{哈佛大学，美国马萨诸塞州剑桥市, MA02138, Arthur\_Jaffe@harvard.edu}、 薛博卿\footnote{上海科技大学, 中国上海市, 201210, xuebq@shanghaitech.edu.cn}}

\date{}

	\maketitle

  {\small \textbf{Abstract:} The ``Millennium Prize Problems'' have a place in the history of mathematics.  Here  we tell some little-known anecdotes from the perspective of the planner of that project. These stories are far from their end; more likely they are just at their beginning.}	
  
	{\small  \textit{摘要：~ “千禧年大奖难题”注定是数学史上经久不息的话题。本文以该事件筹划者的视角讲述一些鲜为人知的趣闻，故事还远没有结束，或许只是刚刚开始。}}
	
	\bigskip
	
	千禧之际，万象更新，数学界的柔风细雨中惊响起初夏的雷鸣。七个重要的数学问题！柒佰万美元的巨额奖金！克雷数学研究所公布的大奖难题在巴黎街头响彻云霄，引人一瞥惊鸿，这一天，正是$2000$年$5$月$24$日。此事缘何而起？它又将何去何从？本文以只言片语，讲述些许这注定不凡的往事。
	
	“千禧年大奖难题”的诞生与克雷数学研究所\footnotemark[1]\footnotetext[1]{Clay Mathematics Institute，它是一个私人非营利机构，于1998年由Landon T. Clay先生资助创立，由Arthur Jaffe教授担任创始会长，致力于促进和传扬“数学思想的美、力量与普适性”(the beauty, power, and universality of mathematical thought)。}（后文简称CMI）的成立密不可分。1998年4月15日的响午时分，哈佛教师活动中心\footnotemark[2]\footnotetext[2]{Harvard Faculty Club, 20 Quincy St, Cambridge, MA 02138.}古朴的餐厅中萦绕着人们的莺莺絮语。此时，Landon Clay先生\footnotemark[3]\footnotetext[3]{Landon Thomas Clay，1926.3.12.-2017.7.29.，美国投资家，1950年毕业于哈佛大学，主修英语文学，他是科学和教育的大力赞助者，克雷数学研究所创始人。}（后文简称LC先生）首次向我\footnotemark[4]\footnotetext[4]{本文中的“我”均指代第一作者Arthur Michael Jaffe，数学物理学家，美国科学院院士，时任美国数学会会长，是克雷数学研究所成立和“千禧年大奖难题”设立的主要筹划人。}
提起了创建一个软件基金会的想法，他部分拥有一家濒临倒闭的公司，正适合改建成基金会。“这是一个争取税收优惠\footnotemark[5]\footnotetext[5]{根据美国法律，缴税时可减免给独立基金会捐赠金额的三分之一。}的好方法，但并不是花钱的最佳途径。”我略作沉思提出了自己的建议，“倘若你愿意为数学作一番事业，我必鼎力相助！”LC先生眼眸一亮，他素来是科学和教育的大力赞助者。大约六周以后，他作出重大决定：另外设立一项基金会，专门支持数学。正所谓兵马未动，粮草先行，这使我备受鼓舞。
	       
	经过深入的思考和充分的准备，我提交了十项可行的项目方案，包括其中的第八项——“千禧年大奖”计划。当时，“千禧年”这个话题在全球如火如荼，我对这一计划青睐有加。而Alain Connes、Andrew Wiles和Edward Witten 这几位享誉世界的数学家们（后文分别简称为AC教授、AW教授和EW教授）也陆续加入了这项事业。同年9月25日，CMI成为注册在美国特拉华州的基金会。第一次董事会议\footnotemark[6]\footnotetext[6]{共有四位参会者：三位董事Clay夫妇和我，以及秘书Barbara Drauschke女士。}在11月10日召开，通过了这十个项目，并选举出科学顾问委员会。
	
	1999年5月10日是看似寻常却又无比奇妙的一天，大约$450$名数学家“不经意间”相聚在麻省理工学院\footnotemark[7]\footnotetext[7]{Massachusetts Institute of Technology.}（后文简称MIT），而CMI的落成典礼也在节日般的气氛中举行了。当天，著名艺术家兼数学家Helaman Ferguson先生\footnotemark[8]\footnotetext[8]{《Helaman Ferguson: Mathematics in Stone and Bronze》（见图a）一书中收录了许多HF先生的杰作，包括多种纽结的雕塑。}（后文简称HF先生）亲手揭开他的杰作——一个阴雕的“$8$字结”雕塑（见图b）——它由花岗岩制成，足足有半吨重。在纽结理论中，“$8$字结”是一个非常重要的范例\footnotemark[9]\footnotetext[9]{Figure eight knot，它是唯一的交叉数为$4$的纽结，并且是素纽结。}，因而它被定为了CMI的标志（见图c）。CMI的主要目标和宗旨为：促进和传播数学知识；在广大科学工作者中宣扬数学领域的新发现；鼓励具备天赋的年轻人从事数学职业；以及对数学研究中的非凡成就或巨大进步进行官方认证。伴随着一系列的公众讲座和相关活动，CMI正式拉开了历史的一角。
		\begin{figure}[h]
		\center
		\subfigure[(a) HF作品集]{
			\includegraphics[width=4.3cm,height=3.2cm]{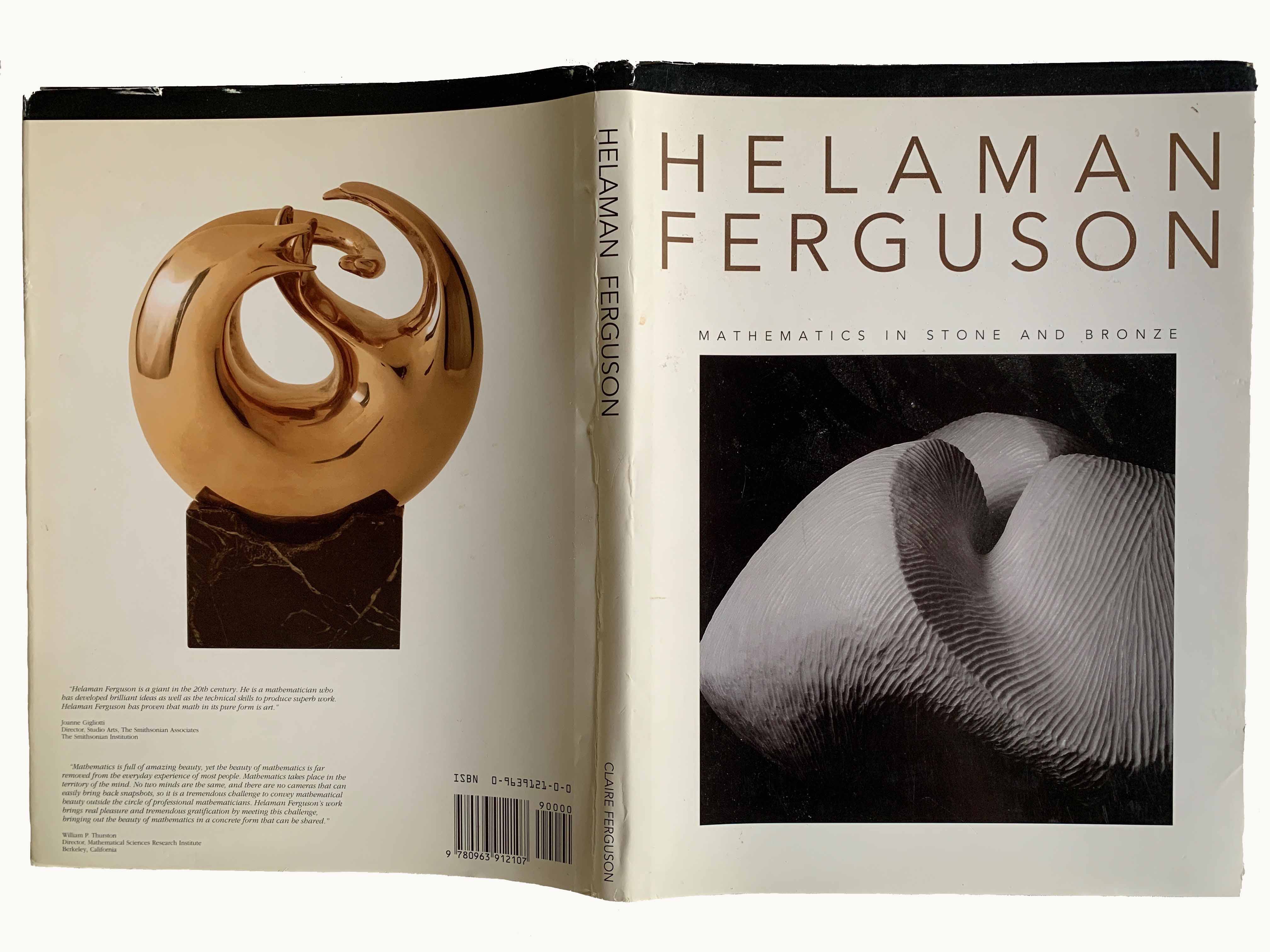}}
		\hspace{0in}
		\subfigure[(b) HF先生和“8字结”阴雕]{
			\includegraphics[width=4cm,height=3cm]{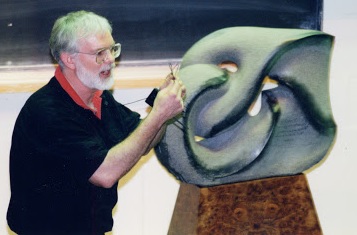}}
		\hspace{0in}
		\subfigure[(c) 如今的CMI标志]{
			\includegraphics[width=3cm,height=3cm]{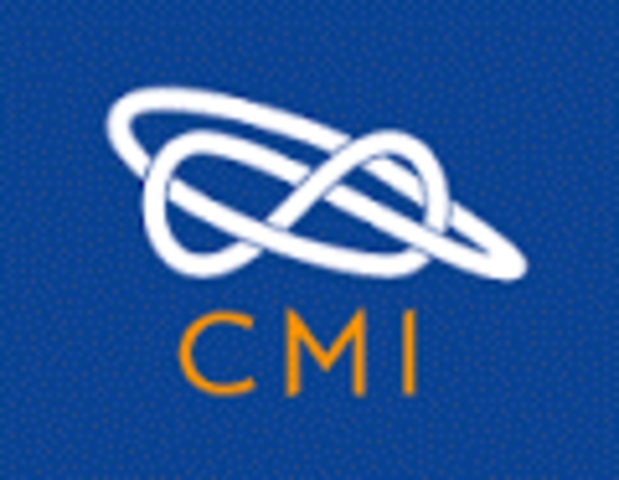}}
	\end{figure}

	“千禧年大奖难题”的诞生还与一百年前珍贵的历史遥相呼应。1900年8月8日，希尔伯特\footnotemark[10]\footnotetext[10]{David Hilbert, 1862年1月23日-1943年2月14日，德国数学家，是十九世纪及二十世纪前期最全面且最具影响力的数学家之一，发明了大量的思想观念（例如不变量理论、公理化几何等），并将哥廷根大学建设为世界数学中心。}在巴黎举行的第二届国际数学家大会上宣布了十个问题，并在1902年正式出版了二十三个问题，它们被合称为“希尔伯特问题”，引领了二十世纪数学的诸多研究方向。需要说明的是，“千禧年大奖难题”的设立初衷与“希尔伯特问题”截然不同，并非是为了预测或影响新世纪数学的发展倾向，而仅仅是为了彰显一些未解的难题，以吸引大众对数学的关注。当然，我们希望把2000年春夏和法国巴黎分别作为这些难题正式公布的时间和地点，以此来表达对希尔伯特这位数学巨匠最诚挚的敬意。

	关于第八项方案的最初计划是通过竞争角逐的方式确立大奖难题，计划书的原文是这样写的：“先择取五十个科研问题，成文并出版成书籍，再从中选取数个‘特别问题’，数量不超过十二个。”一位参与者显得忧心忡忡：“如果公开选取这十二个特殊问题，很容易引发选取流程的政治化倾向并引起争议。”我们最终决定由CMI的科学顾问委员会直接决定问题择取的事宜。该委员会的成员正是前文所提及的AC、AW、EW三位教授\footnotemark[11]\footnotetext[11]{顺序按照英文姓氏的字母序排列。}和我。这时，时光的脚步已然悄无声息地来到了1999年11月，我们猛然发现，距离2000年春夏只剩数月，只有创造奇迹，才可能顺利达成愿景。
	
	难题的数目其实并没有预先设定，我们仅仅认为十二个是一个合理的上限。因为数量足够少，令人得以专注其中；数量亦足够多，使得选题不至于太狭隘。而具体的数目将取决于委员会的工作进程，当时的我们亦无法预测大奖难题的最终面貌。一开始，每个委员会成员都提交了一份清单，所有清单上都赫然出现了两个名词——“黎曼(Riemann)假设”和“庞加莱(Poincar\'{e})猜想”。毋庸置疑，它们入围了。接着，委员会召开了频繁的电话会议，我们兴致勃勃地讨论起各式各样的问题，取得一致意见后，便把一个难题加入到名单中，或者用它替换掉一个已在名单上出现的难题。要知道，讨论这些艰深的数学问题并作出合适的判断是一个充满艰难坎坷的过程，这耗费了我们大量的时间和精力。当遇到一个超出我们专业的问题时，我们会联系相关专家进行咨询。我们努力地让大奖难题不特别集中在某一个数学的领域，也正是由于这个原因，许多同样好的问题遗憾地遭到了淘汰。慢慢地，讨论的效率变得越来越低，再往名单中增加或替换单个难题都变得举步维艰，而我对时间进度也越发担忧。当我们在电话中难以取得进展时，我当机立断，宣布讨论环节结束。这时，刚好列出了七个难题。数字“七”或许有很多特殊的含义，然而这一次，它的出现只因天意。
	
	另一个重要的议题是：这些问题该以什么形式提出？比如，庞加莱猜想可被推广为Thurston猜想，而黎曼假设与Langlands纲领有关。委员会作出决定：难题都\footnotemark[12]\footnotetext[12]{除了后文中的难题(7)以外。}尽可能地以最简单明了的形式呈现。最后，我们还决定了由哪些专家写下难题的具体陈述。这些专家临危受命，在极短的时间内完成了高质量的工作。他们写下的命题以及相关的综述性短文及时地通过了其他专家的审阅、修改及委员会的最终检验。
	
	七个难题如下，顺序依照难题英文名称的字母序排列，括号中所列的是写下难题之具体陈述的作者姓名。
	
	(1) Birch-Swinnerton-Dyer猜想 (Andrew Wiles);
	
	(2) Hodge猜想 (Pierre Deligne);
	
	(3) Navier-Stokes方程解的存在性与光滑性 (Charles Fefferman);
	
	(4) P/NP问题 (Stephen Cook);
	
	(5) 庞加莱猜想 (John Milnor);
	
	(6) 黎曼假设 (Enrico Bombieri);
	
	(7) Yang-Mills规范场存在性与质量间隙 (Arthur Jaffe和Edward Witten)。

	这七个难题中，黎曼假设是唯一一个在“希尔伯特问题”里就出现过的，已历经一百多年却仍然巍然屹立。它在山顶的风景如此令人迷醉，山间的雾霭却乱人眼眸，人们手执大斧披荆斩棘，却始终无法找到一条通往山巅的清晰途径。
	
	需要特别强调的是，这是七个（当时）未解决的重要问题，但并非“最”重要。无论是揭开古老的未解之谜，还是发现全新的研究方向或领域，都无比艰难。由前者所取得的成果较易为当世之人（尤其是数学家们）所敬仰，而后者的成就往往需要经过更多的时间积淀才能被世人所理解和接受，两者都难能可贵。

	关于获奖的规则细节，以下几条是原则性的：获奖者必须证明或证否其中一个难题；解答不能直接提交给CMI，需在正式的学术期刊上发表；CMI将在解答被发表两年之后启动审查程序；对于合作取得的成果，或是重要的先验想法，时任董事会将根据科学顾问委员会的建议来决定\footnotemark[13]\footnotetext[13]{在当时，第一作者是董事会中唯一的数学家。}他们应分享的荣耀。
	
	奖金金额的设定过程则是另一个非常有趣的故事。我的原始想法是设立滚动奖池制：基金会每年向奖池中投入一定量的金额，当一个大奖难题被解决时，获奖者们将分享的奖金为
	\[
	\text{奖金} \,\, = \,\,\frac{\text{当时的奖池总额}}{\text{当时剩余的未解决难题数量}+1}.
	\]
	按照这样的设置，如果一个难题很快就被攻破，它产生的奖金是相对朴素的。如果一些难题一直悬而未决，奖池将会越滚越大，最后变成超级巨奖。恰好在2000年春的一份《泰晤士报》\footnotemark[14]\footnotetext[14]{Times of London.}中，英国费伯出版社\footnotemark[15]\footnotetext[15]{Faber\&Faber.}为了提高一本书的关注度，为解答哥德巴赫猜想\footnotemark[16]\footnotetext[16]{Goldbach's conjecture，数论中存在最久的未解问题之一，其表述为：“任一大于$2$的偶数可以表示成两个素数之和。”}提供了一百万美元的奖金。在当时，百万美元可是相当大的一笔财富。然而，出版商只给这笔大奖设置了两年的有效期，他们并没有预备那么多钱，而是就此奖金向劳合社\footnotemark[17]\footnotetext[17]{Lloyd's of London.}购买了保险。对一个古老的数学难题而言，两年的时间如过眼云烟，该出版商用极低的风险换得了人们的极大关注。在这样的背景下，我提议为这七个大奖难题预设七百万美元的总奖金。略微遗憾的是，奖池通过投资逐年递增的提议被CMI的律师驳回了。受限于筹备时间的仓促，我们对奖池问题没能进一步地协商，这就是七百万美元大奖的由来。

	在甄选难题和设置奖金的同时，另一项艰巨的任务是筹备将在巴黎召开的千禧年大会。AC教授作为法国方面的联络人，联系了法兰西公学院\footnotemark[18]\footnotetext[18]{Coll\`{e}ge de France.}作为大会的东道主。大会被定于2000年5月24-25日召开，大奖难题的发布安排在24日下午，而25日则是学术报告会，报告人均为与CMI相关的青年学者，包括Terrence Tao、Manjul Bhargava、Dennis Gaitsgory等，彼时他们方弱冠少年，未及而立之年。另外，出于保密和出其不意的目的，我们给大会起了一个波澜不惊的名字——CMI第二次年度会议。一般而言，举办一个大型国际会议需要提前相当多的时间，以便参与者们提前做好行程规划。由于筹备开始的时间太晚，我们无可避免地陷入了一些困境：在2000年春夏之间，无论时间怎么调整，都会与一些其他的数学会议相冲突。比如，二月份时我们曾希望邀请Raoul Bott和Jean-Pierre Serre两位教授来报告这七个大奖难题，而他们被授予了沃尔夫数学奖\footnotemark[19]\footnotetext[19]{Wolf Prize, 由Ricardo Wolf创立的沃尔夫基金会颁发，奖励对推动人类科学与艺术文明做出杰出贡献的人士，每年评选一次。其中的数学奖被认为是“终身成就奖”。}，颁奖典礼与千禧年大会的时间是重叠的。另外，由于无法预估参会人员的规模，在大会召开前几天，AC教授决定为每位注册参会者签发门票，每张门票上都印有专门的座位号。对于届时未能进场的观众，我们安排了分会场收看直播。AC教授以及大洋两岸诸多赋有奉献精神的支持者们，投入了极大的工作热情，最终我们克服了所有困难，保障了“千禧年大奖”项目的成功实施。

	法兰西公学院首脑Gilbert Dagron教授（下文简称GD教授）相信“千禧年大奖难题”将会是科学史上的一项重大事件，他为我们提供了不可或缺的支持（见图l中信件）。更为重要的是，我们被允许借用他杰出的助手V\'{e}ronique Lema\^{\i}tre女士（下文简称VL女士）。这位令人惊喜的VL女士与巴黎地区的科学记者们保持着密切而友好的关系，她为我们安排了各类采访。AC教授在大会前一周接受了《世界报》\footnotemark[20]\footnotetext[20]{Le Monde.}的采访，而我也提前三天接受了美联社\footnotemark[21]\footnotetext[21]{Associated Press.}的访谈。在难题发布的当天早晨，VL女士还召集了媒体发布会，吸引了三十多位科学记者参加。自然而然地，所有的消息都被设置了“禁令”，以防在官方发布之前就遭到曝光。

	法兰西公学院的玛格丽特·德·那瓦尔圆形剧场\footnotemark[22]\footnotetext[22]{Amphith\'{e}\^{a}tre Marguerite de Navarre, Coll\`{e}ge de France.}刚落成不久，宽大的舞台后设置了两个大型的投影巨幕，环绕着的观众座椅则层叠而上，如孔雀开屏般展成一个扇形。5月24日下午2时，全场座无虚席。大会第一项议程是克雷研究奖\footnotemark[23]\footnotetext[23]{Clay Research Award.}颁奖仪式，获奖人为Laurent Lafforgue教授（下文简称LL教授），Lavinia Clay夫人为他颁发了特制的奖杯——一个小型的“8字结”青铜艺术品（见图e）。接着，在观众们惊奇的眼神中，我从讲台后方取出了藏着的第二个奖杯：“正如数学中的对称性、多重性！我们有两个奖杯！”第二个奖项被授予了事先并不知情的AC教授。LL教授因为有关Langlands纲领\footnotemark[24]\footnotetext[24]{Langlands program，这是数学中一系列影响深远的构想之一，联系了数论、代数几何与约化群表示理论，一些数学家将其描述为数学中的“大一统”理论。}的工作获奖，其背后还伴随着一些轶事。大会之后没多久，他给我写信说想退回这一奖项，因为当他在准备亨利-庞加莱研究所\footnotemark[25]\footnotetext[25]{Institut Henri Poincar\'{e}.}演讲的过程中，发现了证明中的一个严重漏洞，因而惶恐不安。我们都给他以鼓励，并让AW教授安慰他，因为AW教授曾有过类似经历。正如我们所料，不久后LL教授便修补了漏洞，得到了完美的证明。

	\begin{figure}[h]
		\center
		\subfigure[(d) 大会海报]{
		\includegraphics[width=2.5cm,height=3.5cm]{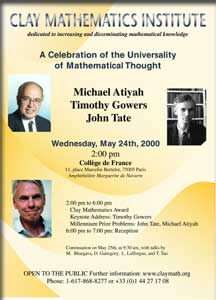}}
		\subfigure[(e) 2000年的CMI标志兼克雷研究奖奖杯]{
		\includegraphics[width=2.5cm,height=2.5cm]{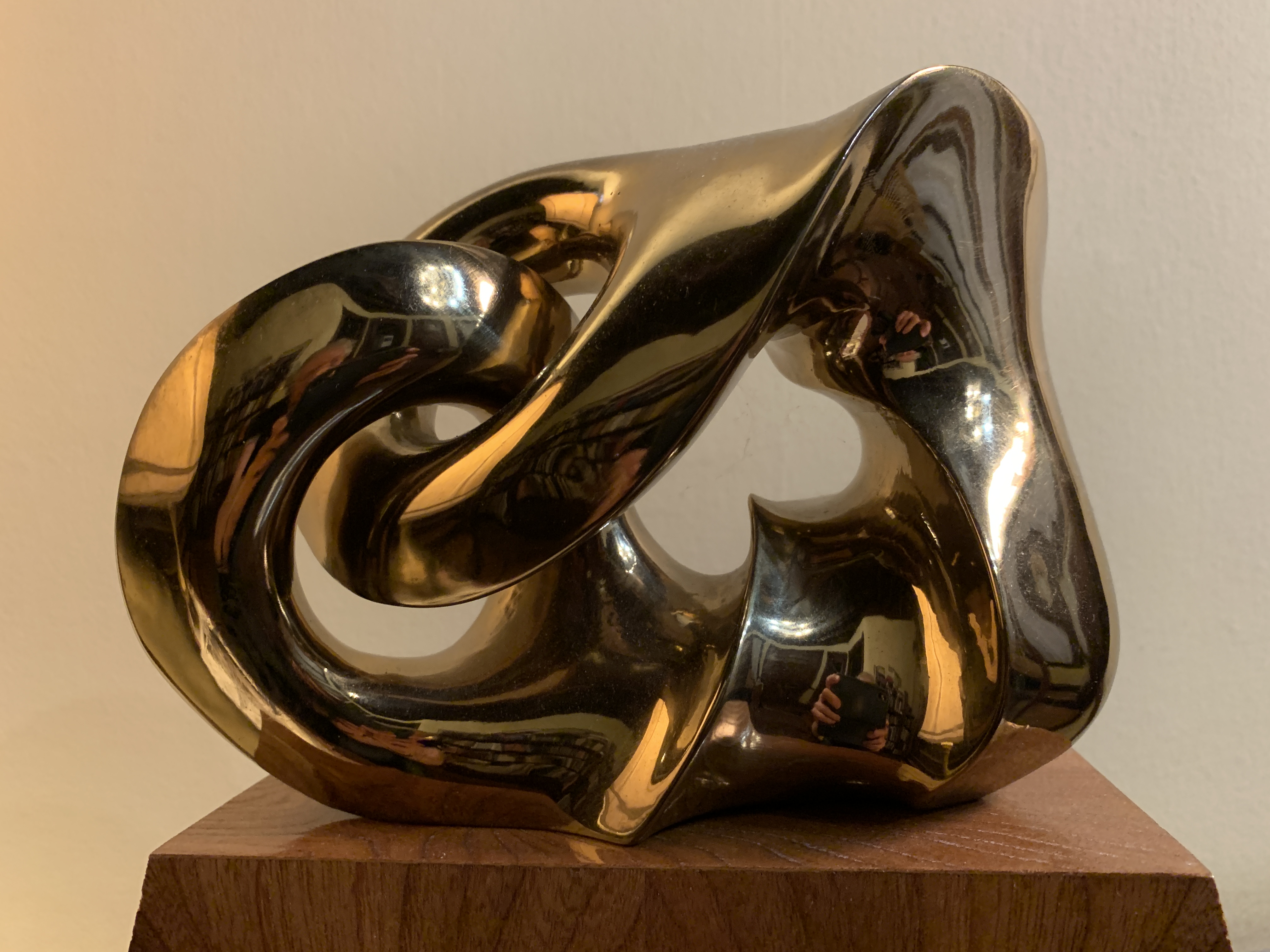}}
		\hspace{-0.1in}
		\subfigure[(f) 大会现场]{
		\includegraphics[width=3.5cm,height=2.5cm]{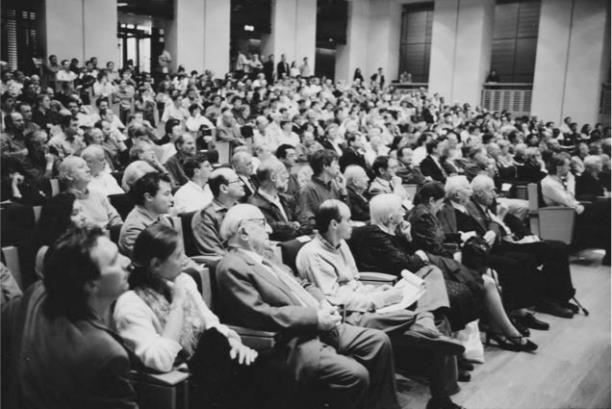}}
		\hspace{-0.1in}
		\hspace{-0.1in}
		\subfigure[(g) 录像集]{
		\includegraphics[width=2.5cm,height=2.5cm]{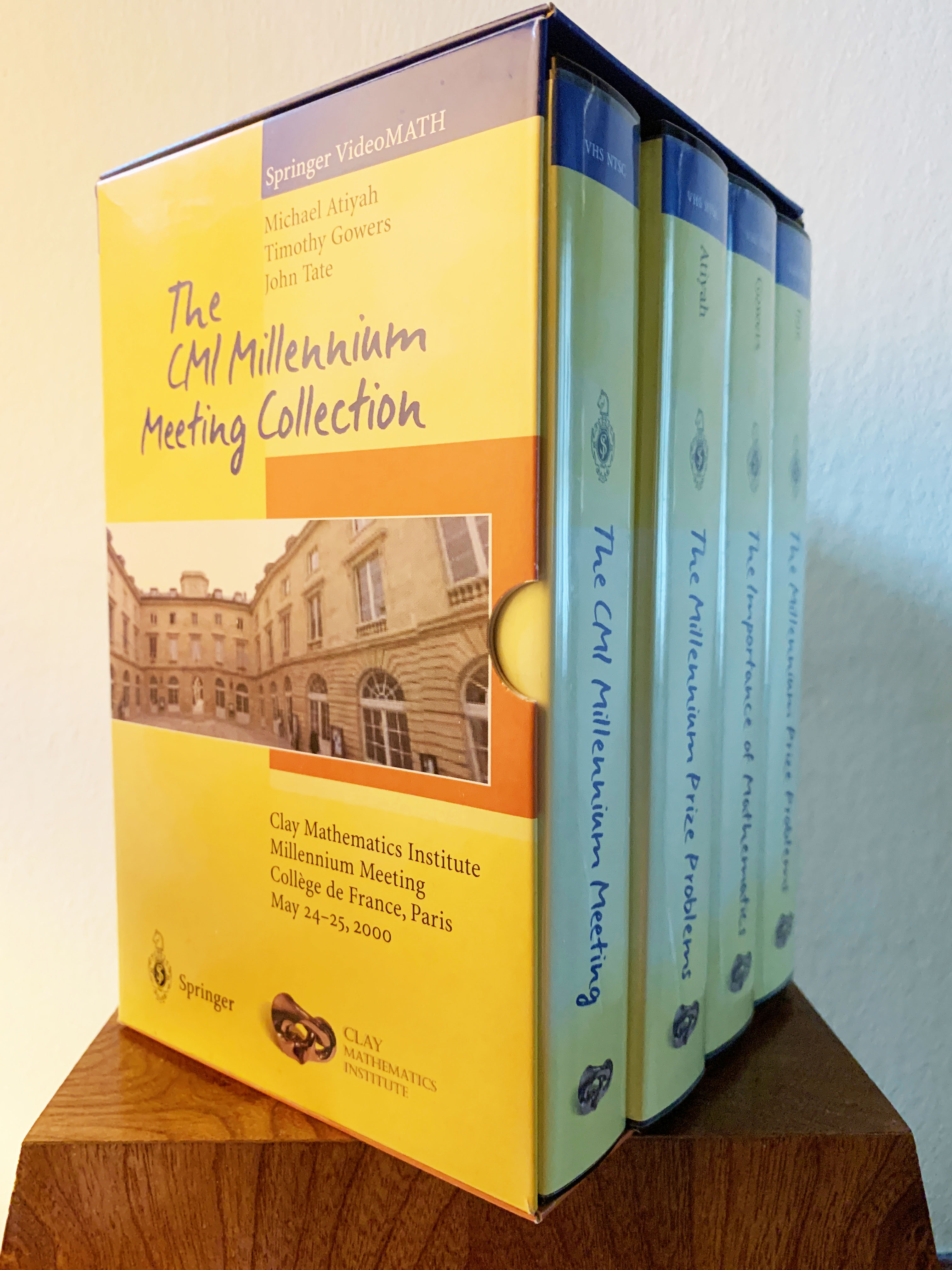}}
	\end{figure}
	
	颁奖之后，Timothy Gowers教授作了题为“数学如此重要”的主题演讲\footnotemark[26]\footnotetext[26]{Keynote Address: The Importance of Mathematics.}，而Michael Atiyah和John Tate两位教授则逐个介绍了这七个难题（见图d中海报）。那些令人一头雾水的专业术语被酝酿成了极为朴实的言语，而难题的来龙去脉如行云流水一般，成了他们口中淡淡倾吐着的动人故事。午后的时光转瞬即逝，观众们聚精会神，听得如痴如醉（见图f）。随着现场播放的录音中响起希尔伯特低沉而充满张力的德语，全场报以热烈的掌声，这是他在1930年发表的退休演说。至此，尚没有任何人提起过攻破这七个难题是有奖金的。到最后，我邀请大家享用招待会的香槟和糕点，并轻描淡写地说：“顺便提一句，如果今晚你恰好解答了其中的一个难题，请向期刊投稿，我们将在两年后颁奖，每个难题的奖金是：一百——万——美元！”刹那间，惊呼声、口哨声、击掌声，震天撼地！议论声纷呈而起，久久不曾散去。
	
\begin{window}[0,l,{\mbox{\includegraphics[width=2.5cm,height=3cm]{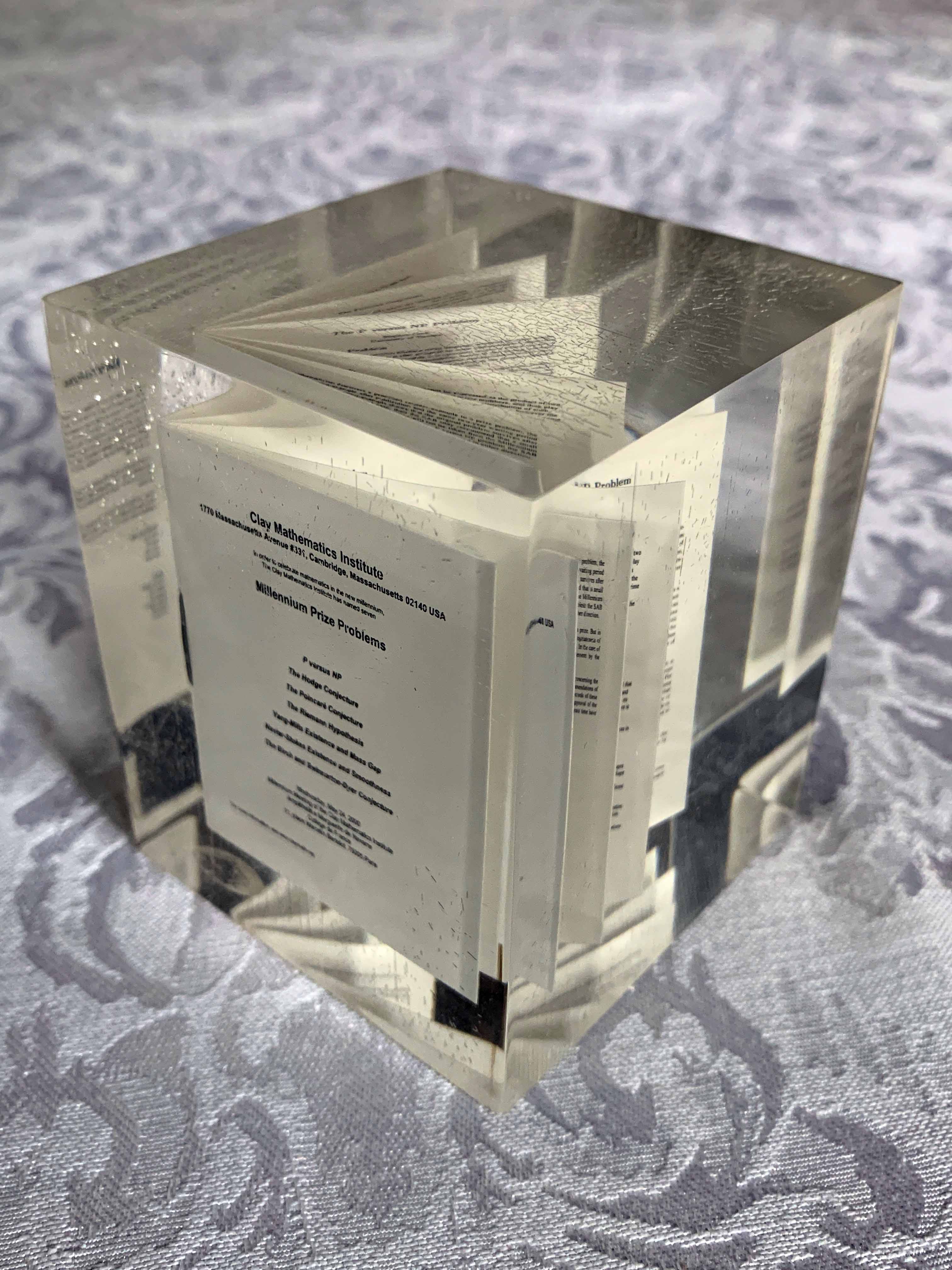}}}, {\footnotesize (h) 立方体纪念品}]
	大会期间还有一些插曲，包括铜管乐团优美的演奏，以及LC先生、GD教授和一名政府要员的讲话，而大会之后安排了招待会和晚宴\footnotemark[27]。在晚宴上，激动的HF先生向大家阐释着有关“8字结雕塑”的灵感源泉。“这是一个鞍点。”他指着右手的虎口，旋转着的手掌划起蜜蜂的$8$字舞，“看！就如这般，多么神奇的结构啊！”每位晚宴的客人都收到一枚晶莹剔透的立方体艺术品（见图h），边长6.5厘米，它似琥珀一般封存了一本“千禧年大奖难题”的小册子，如此典雅精致的小礼物赋予了这个巴黎之夜更深的诗情画意与惊喜。
\end{window}
\footnotetext[27]{请见关于千禧年大会的录像集《The CMI Millennium Meeting Collection》（见图g），共四个录像，由Fran\c{c}ois Tisseyre 制作，Springer 出版社出版，ISBN：978-3-540-92657-3。这些录像也可以在网上找到。}

\begin{window}[2,r,{\includegraphics[width=2.5cm,height=3.5cm]{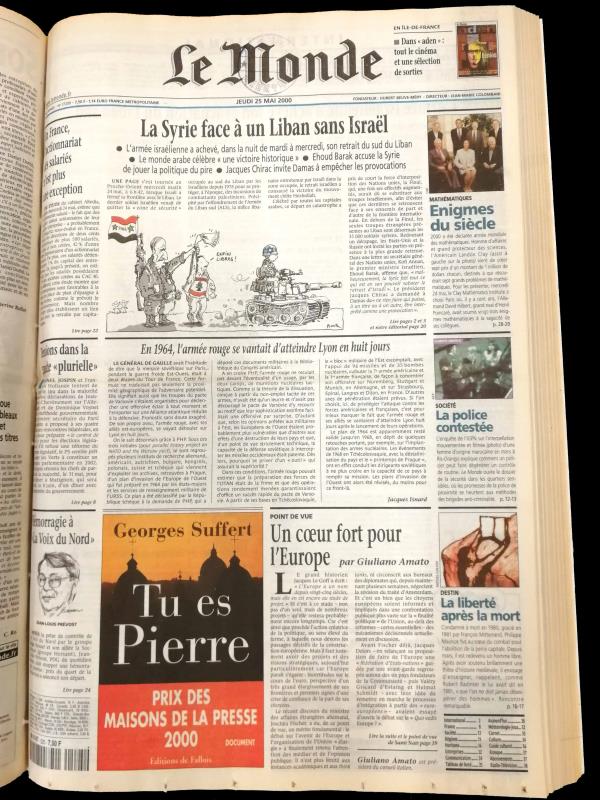},\includegraphics[width=2.5cm,height=3.5cm]{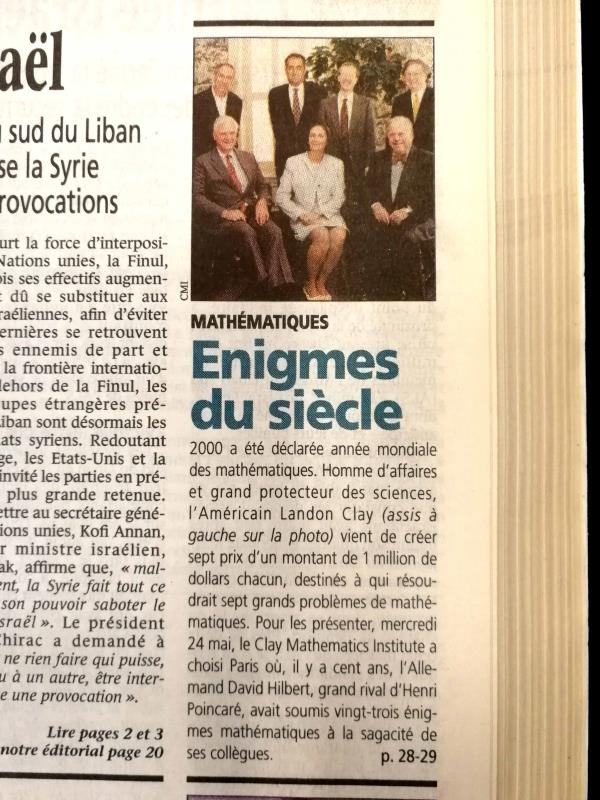}}
, {\footnotesize ~~ (i) 2000年5月25日《世界报》首页}]
“千禧年大奖难题”公布后的即时反应异常热烈。招待会结束后，我们乘坐大巴车前往晚宴，在车上，VL女士给我们分发了《世界报》刚刚发行的报纸。该报社每天傍晚刊发次日的报纸，因此日期显示为5月25日。在这份报纸中，科学记者Jean-Fran\c{c}ois Augereau写下了四篇关于“千禧年大奖难题”的不同文章。其中，首页右上方（见图i）刊登了一张包含CMI董事会和科学顾问委员会的照片（见图j），这是报社自己在网上找到的，摄于CMI落成典礼的当天。由美联社的科学记者Jocelyn Gecker撰写的长文则刊发在25日世界各地的多家报纸上。这些报道引申出后续成百上千的故事，许多电视、网络媒体都争相跟进。
\end{window}
	
	\begin{figure}[h]
		\center
		\subfigure[(j) CMI董事会与科学顾问委员会]{
			\includegraphics[width=9cm,height=6cm]{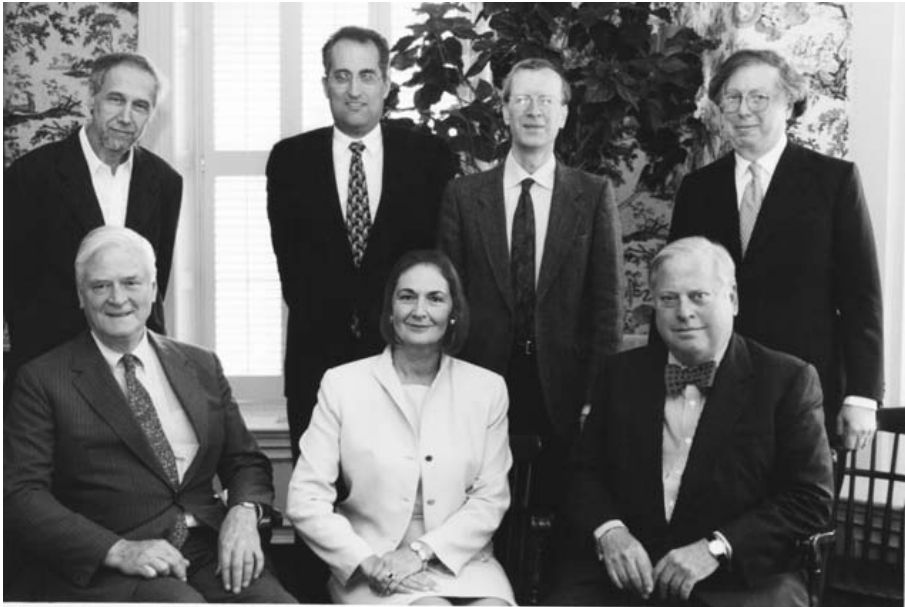}}
			\vspace{-0.5cm}
		\caption*{ \footnotesize   \indent ~~~~~~~~~~~~~~第二排左起：Alain Connes, Edward Witten, Andrew Wiles, Arthur Jaffe\protect\\ \indent ~~~~~~~~~~~~~~~~第一排左起：Landon T. Clay, Lavinia D. Clay, Finn M.W. Caspersen
		}
	\end{figure}
	
	哈佛大学本科生Kelvin Liao当时在CMI做兼职，凭借着电脑技术和艺术两方面都极高的天赋，他设计了CMI网站的前三个版本。当大会在法国巴黎公布“千禧年大奖难题”的同一时刻，他在美国马萨诸塞州上线了CMI网页的第二版，同时贴出了“千禧年大奖难题”的官方公告。公众的反响远远超出了我们的预期，公告刚一发布，CMI网站便因为过载的流量而崩溃了。我联系了美国数学会（后文简称AMS）的执行董事John Ewing先生（后文简称JE先生）,他决定在AMS的服务器上建立CMI网站的镜像。要知道AMS网站除了发布数学新闻之外，还承载了大量的电子书籍和期刊数据，并为全世界范围的各类数学事务提供服务，它的容量和带宽远胜CMI所使用的商用服务器。两天后我返回美国，JE先生立马给我打来电话，他忧心忡忡：“AMS网站要崩溃了！我们绝不能失去AMS的杂志和书店！”“能否再坚持一天？”我也万般无奈，“可能马上就趋于稳定了？”幸运的是，AMS的服务器顶住了冲击，经过大约一周的时间，网络点击量趋于平稳。后来，我通过MIT的工程学教授Charles Leiserson联系了阿卡迈科技公司\footnotemark[28]\footnotetext[28]{Akamai Technologies.}，拯救了这一危机。
	
	经过成千上万次报纸、杂志、电视的曝光，数学收到了各界人士的广泛关注。比如在一档晨间电视台的优质栏目中，ABC电视台的新闻台柱之一Chuck Gibson和MIT的Michael Sipser教授一起给大家讲述了P/NP问题。有许多年轻的学生说他们对数学倍感兴趣，并且想了解这些大奖难题。而“千禧年大奖难题”产生的长远影响之一是：大量相关的数学书籍问世，其中涉及黎曼假设和庞加莱猜想的举不胜举。一些人着迷于尝试解决其中一个难题，然而赢得百万美元的奖金并不容易。另外，正如之前所说，解答著名的古老问题只是数学研究的一个方面，探究未知数学的新领域或探寻未来数学的新视角也同等重要。

\footnotetext[29]{arXiv.org，这是一个收集物理学、数学、计算机科学、生物学等领域的论文预印本的网站，提供免费和开源的服务。它不是学术期刊，在其上张贴论文预印本无需经过同行评审。}
\footnotetext[30]{The State University of New York at Stony Brook.}
\footnotetext[31]{Fields Medal，以加拿大数学家John Charles Fields命名，又名国际杰出数学发现奖，在四年一次的国际数学家大会上颁发给2-4名有卓越贡献且年龄不超过40 岁的数学家。}
\footnotetext[32]{Steklov Mathematical Institute.}

\begin{window}[3,r,{\mbox{\includegraphics[width=4cm,height=3cm]{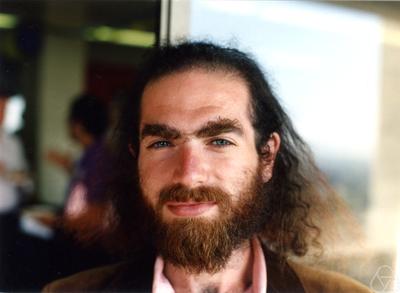}}}, {\footnotesize ~~~~~~~~~ (k) GP博士}]
七个大奖难题之一——庞加莱猜想——如今已被完全解决，这是一件在数学界耳熟能详的传说。在2002和2003年，Grigori Perelman博士（见图k，下文简称GP博士）在预印本网站\footnotemark[29]张贴了三篇论文宣称解决了庞加莱猜想的推广形式——Thurston猜想。2003年在石溪大学\footnotemark[30]，我曾与他聊了一整天，这一过程令人陶醉。我告诉了他CMI奖金的事，而他对此仅有的评论是：“金钱很危险！”经过数个数学研究团队的补充和验证，到2006年时数学界已认可GP博士的工作无显著漏洞。这时GP博士不仅放弃了百万美元的奖金，还拒绝了包括菲尔兹奖\footnotemark[31]在内的所有奖项。他甚至辞去了斯捷克洛夫数学研究所\footnotemark[32]的工作，与母亲一起过着贫困的生活。在2008年，我再次尝试联系他，并没有成功。听说他最大的人生乐趣是听音乐，并与音乐家们交往。
\end{window}

一颗耀眼的明星就这样藏匿起了自己的光芒，大隐隐于寰宇之间。而下一段传奇的主人公，会不会是正在阅读此文的你呢？故事还远没有结束，或许只是刚刚开始。

	\begin{figure}[H]
		\center
		\subfigure[(l) GD教授的信件]{
		\includegraphics[width=9cm,height=12cm]{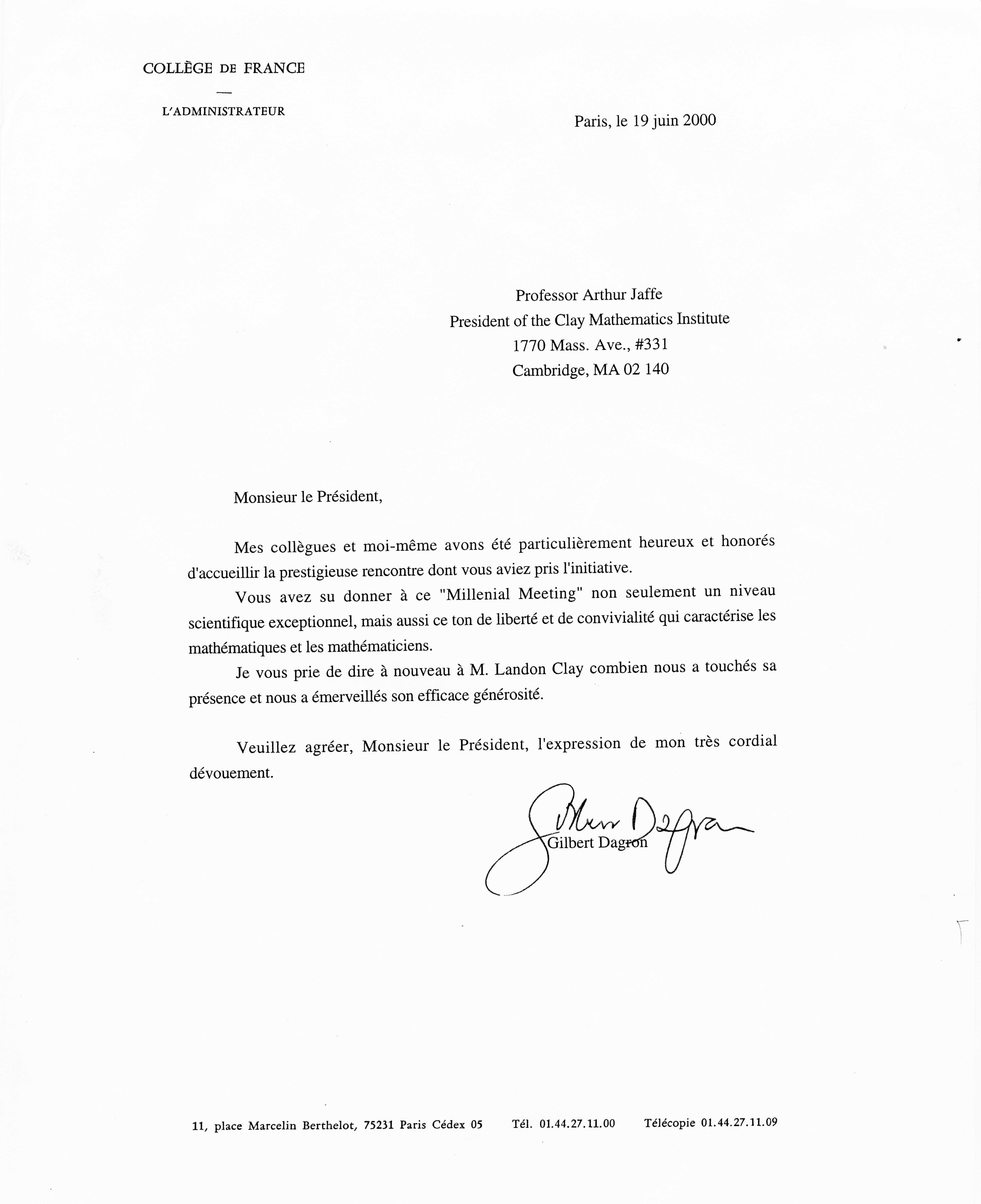}}
	\end{figure}

	本文所叙述的往事以第一作者的演讲\footnotemark[33]\footnotetext[33]{时间：2018年7月27日10:00；地点：中国科学院数学与系统科学研究院N219；报告人：Arthur Jaffe院士；报告题目：Story Behind the Millennium Prize Problems in Mathematics.}和纪实报告\footnotemark[34]\footnotetext[34]{A. Jaffe, The Millennium Grand Challenge in Mathematics, Notices of the American Mathematical
		Society, 2006, 6(53), 652-660.}为基础，由第二作者进行了少许演绎，基本都保持了事件的原貌。文中的第一人称“我”指代第一作者。所有当代人物的姓名均以其原始的语言展示，其余的人名、地名、机构名称等均译成了中文。
	
	\medskip
	
	{\footnotesize 致谢：由衷感谢任云翔博士，他为修缮及翻译本文提供了不可或缺的帮助。亦特别感谢张睿燕女士，她帮助作者们寻找到有关《世界报》的重要史料。本文由Templeton Religion Trust(TRT 0159)和国家自然科学基金(No. 11701549)部分资助。}


\clearpage\end{CJK*}
\end{document}